\documentclass[reqno]{amsart} \usepackage{amscd, amssymb, epsfig}
\newtheorem{theorem}{Theorem}[section]
\newtheorem{proposition}[theorem]{Proposition}
\newtheorem{lemma}[theorem]{Lemma}
\newtheorem{corollary}[theorem]{Corollary}
\theoremstyle{definition}
\newtheorem{definition}[theorem]{Definition}

\theoremstyle{remark} \newtheorem{remark}[theorem]{Remark}
\numberwithin{equation}{section}
\newcommand{\field}[1]{\ensuremath{\mathbb{#1}}}
\newcommand{\CC}{\field{C}}

\newcommand{\RR}{\field{R}}

\newcommand{\diff}[1]{\mathcal{#1}} 
\DeclareMathOperator{\id}{id} 
\DeclareMathOperator{\tr}{tr}

\DeclareMathOperator{\im}{Im} 
\DeclareMathOperator{\Aut}{Aut}
\newcommand{\del}{\partial}
\newcommand{\delb}{\bar\partial}

\newcommand{\ot}{\otimes}

\newcommand{\e}{\varepsilon}
\newcommand{\g}{\gamma}
\newcommand{\G}{\Gamma}

\newcommand{\p}{\partial}

\begin{document}
\begin{flushright}
\end{flushright}
\title[Deformation Quantization]
{Deformation Quantization of K\"{a}hler Manifolds}
\author{Nikolai Reshetikhin} \address{Department of
Mathematics\\ University of California at Berkeley\\Berkeley, CA 94720 \\ USA}
\email{reshetik@math.berkeley.edu}
\author{Leon A. Takhtajan} \address{Department of Mathematics \\
SUNY at stony Brook \\ Stony Brook, NY 11794-3651 \\ USA}
\email{leontak@math.sunysb.edu}
\begin{abstract}

We present an explicit formula for the deformation quantization on K\"{a}hler
manifolds.
\end{abstract}
\maketitle
\begin{quote}
{\it Dedicated to our teacher Professor L.D.~Faddeev on his 65th birthday}
\end{quote}
\section{Introduction}
The concept of quantization has a long history. Mathematically, it originated
from Hermann Weyl~\cite{weyl}, who introduced a mapping from classical
observables --- functions on the phase space $\RR^{2n}$ --- to quantum
observables --- operators in the Hilbert space $L^2(\RR^n)$. The idea was to
express  functions on $\RR^{2n}$ as Fourier transforms and to replace
the exponential kernel --- the one-dimensional representation of the
Heisenberg group --- by the corresponding unitary group element
in the von Neumann infinite-dimensional representation of the Heisenberg group,
parameterized by the Planck constant $\hbar$. The inverse mapping was
constructed by E.~Winger~\cite{wigner} by interpreting classical observables
as symbols of operators. Later on, Moyal~\cite{moyal} and Groenewold~\cite{gro}
have shown that the symbol of the commutator or of the product of operators is the sine
(respectively, exponential) function of the bidifferential operator of the Poisson bracket
applied to the corresponding symbols.

In the early 70's, Berezin was developing the general mathematical
definition of quantization as a kind of a functor from the category of
classical mechanics to a certain category of associative algebras~\cite{berezin}.
In particular, he considered quantization on K\"{a}hler manifolds. However,
Berezin's approach was analytic in nature and it required very restrictive
geometric conditions on the underlying K\"{a}hler structure~\cite{berezin}.
As a result, the only meaningful examples established in~\cite{berezin} were
$\CC^n$ and homogeneous bounded domains in $\CC^n$. Later on, Berezin's
approach and its relation with geometric quantization were developed
in~\cite{rcg1, rcg2}. Only recently in~\cite{kar1} A.~Karabegov, using local
operator-type constructions, has proved that on an arbitrary K\"{a}hler manifold
there exists a special deformation quantization (with ``separation of
variables'') which is local and is parameterized by an arbitrary formal
deformation of the K\"{a}hler potential.

About the same time as Berezin, Bayen, Flato, Fronsdal, Lichnerowicz and
Sternheimer~\cite{BFFLS, BFFLS2} considered quantization as a deformation
of the usual commutative product of classical observables into a
noncommutative associative product $\star_{\hbar}$ which is parameterized by the
Planck constant $\hbar$ and satisfies  the correspondence principle.
In~\cite{BFFLS,BFFLS2} they systematically developed the concept of
deformation quantization as a theory of $\star$ - products
and gave an independent formulation of quantum mechanics based on this notion.
The BFFLS approach was inspired by the work of Vey~\cite{vey} on the
deformation of Poisson structures, where the idea of the $\star$ - product had
already appeared. Subsequently, the work~\cite{BFFLS,BFFLS2} laid a foundation
for the deformation quantization, a theory that utilizes pure algebraic
approach based on Hochschild cohomology to the deformation of associative
algebras, developed by Gerstenhaber~\cite{ger}, together with the
Chevalley-Eilenberg cohomology of the Poisson-Lie algebra of classical
observables.

Because of its physical origin and motivation, the problem of deformation
quantization was first considered for symplectic manifolds. First results in
this direction were theorems by Vey~\cite{vey} and Neroslavsky and
Vlasov~\cite{nero} that deformations of the Poisson-Lie algebra and of the
commutative algebra of classical observables exist whenever the third Betti
number of a symplectic manifold vanishes. It was shown
in~\cite{lichn} that for any symplectic manifold there exists a fourth
order quantization. Finally, the problem of existence and classification of
$\star$ - products on arbitrary symplectic manifolds was solved by De
Wilde-Lecomte~\cite{dwlc,dwlc-1}. Another approach was taken
in~\cite{maeda, maeda2} and Fedosov~\cite{fed}. In particular,
in~\cite{fed} a rather explicit construction of
the $\star$ - product was given.

The problem of deformation quantization is naturally formulated for the
Poisson manifolds as well. Until recently, the best known example of such
quantization was the Etinhof-Kazhdan theorem~\cite{ek} that every
Poisson-Lie group can be quantized. The situation changed when
M.~Kontsevich proved~\cite{kon} his famous formality conjecture, which,
in particular, implies that deformation quantization exists on
every Poisson manifold. In addition to this general result, Kontsevich
presented a universal formula for the deformation quantization of an
arbitrary Poisson structure on $\RR^n$. This remarkable formula gives a
$\star$ - product of smooth functions on $\RR^n$ as formal power series in
$\hbar$ of bidifferential operators parameterized by graphs, with coefficients
given by integrals over certain configuration spaces. As was
stated in~\cite{kon}, this series resembles a perturbation expansion of the
topological sigma-model of B type (after a suitable supersymmetry
transformation). Recently, this perturbation expansion was carried through
in~\cite{jovanni}.

The present paper appeared from an attempt to understand further the
analogy with the perturbation expansion. Specifically, we observe that
there exists a simple integral formula for the deformation quantization on
K\"{a}hler manifolds based on the finite-dimensional formal integral,
defined through the Laplace method approximation.  Quite naturally,
coefficients of the corresponding formal power series can be written
in terms of ``Feynman graphs".

Remarkably, the integral that we consider is almost identical to the one
used by Berezin~\cite{berezin}. The difference, however, is that now we
regard it as a formal power series in the deformation parameter $\hbar$,
and our approach works for an arbitrary K\"{a}hler manifold, whereas
Berezin was considering a $\star$ - product given by a convergent integral
in some complex domain. As the result, he was forced to impose certain global
restrictions on the K\"{a}hler manifold, which make his approach
meaningful only for $\CC^n$ with the standard K\"{a}hler structure and for
flag varieties, as discussed in~\cite{rcg1,rcg2} (and in~\cite{kar} for
generalized flag varieties). In comparison, our approach gives a direct proof
of the theorem proved by Karabegov~\cite{kar1} in terms of Berezin's original
integral.

The organization of the paper is as follows. In section~ 2 we
recall some basic facts from Berezin's paper~\cite{berezin}. In the next
section we construct the non-normalized $\star$ - product: deformed product
for K\"{a}hler manifolds, given by a simple (formal) integral formula, which
does not preserve the unit element of the commutative algebra of classical
observables. In section~4 we present the normalized $\star$ - product,
which preserves the unit element, and show that there exists formal Bergman
kernel. In the appendix, we prove the equality between two types
of the formal integrals introduced in the paper.

\section{Berezin's quantization}
\subsection{Covariant symbols} \label{formula}

Consider a complex manifold $\CC^n$ endowed with a K\"{a}hler metric
\begin{equation}
ds^2=\sum_{i,j=1}^n h_{i\bar{j}}dz^i\otimes d\bar{z}^j.
\end{equation}
Since $\CC^n$ is a two-connected Stein manifold, every K\"{a}hler metric on it
admits a global potential $\Phi(z,\bar{z})$, i.e.
\begin{equation} \label{Phi}
h_{i\bar{j}}=\frac{\del^2\Phi}{\del z^i\del\bar{z}^j},
\end{equation}
where we are using complex coordinates $z^i, \bar{z}^i,~i=1, \dots, n$, on
$\CC^n$.
Let
\begin{equation} \label{omega}
\omega=-\frac{1}{2}\im ds^2=\frac{\sqrt{-1}}{2}\sum_{i,j=1}^n
h_{i\bar{j}}dz^i\wedge
d\bar{z}^j,
\end{equation}
be the corresponding symplectic form. It induces a Poisson structure on the
algebra $\diff{A}=C^{\infty}(\CC^n)$ of classical observables
\begin{equation}
\{f_1,f_2\}=\frac{2}{\sqrt{-1}}\sum_{i,j=1}^n h^{i\bar{j}}(\del_i f_1\delb_j
f_2-\delb_j f_1 \del_i f_2),
\end{equation}
where $\del_i=\del/\del z^i,~\delb_i=\del/\del\bar{z}^i$ and $h^{i\bar{j}}$
are matrix elements of the inverse matrix $H^{-1}$ of the matrix
$H=(h_{i\bar j})$ of the K\"{a}hler metric.

Let $\hbar$ be a positive real number. Denote by $d\mu_{\hbar}=(\omega/\pi
\hbar)^n$ normalized volume form on $\CC^n$,
\begin{equation}
d\mu_{\hbar}(z,\bar{z})=\det H\prod_{i=1}^n \frac{|dz^i\wedge d\bar{z}^i|}
{2\pi\hbar},
\end{equation}
and define $\diff{H}_{\hbar}$ as the Hilbert space of complex valued, measurable
functions on $\CC^n$ which are square summable with respect to the measure
$e^{-\Phi(z,\bar{z})/\hbar}d\mu_{\hbar}(z,\bar{z})$. The inner product of
$f_1,f_2\in\diff{H}_{\hbar}$ is given by
\begin{equation} \label{sc-prod}
(f_1,f_2)=\int_{\CC^n}f_1(z,\bar{z})\overline{f_2(z,\bar{z})}
e^{-\Phi(z,\bar{z})/\hbar} d\mu_{\hbar}(z,\bar{z}),
\end{equation}
where in order to distinguish $C^{\infty}$ functions on $\CC^n$ from the
holomorphic ones, we denote the former by $f(z,\bar{z})$, and the latter
by $f(z)$.

Let $\diff{F}_{\hbar}\subset\diff{H}_{\hbar}$ be the subspace of holomorphic
functions and $P: \diff{H}_{\hbar}\rightarrow\diff{F}_{\hbar}$ be the
corresponding orthogonal projection operator. Explicitly,
\begin{equation} \label{projector}
(Pf)(z)=\int_{\CC^n}B_{\hbar}(z,\bar{v})f(v,\bar{v})e^{-\Phi(v,\bar{v})/\hbar}
d\mu_{\hbar}(v,\bar{v}),
\end{equation}
where $B_{\hbar}(z,\bar{v})$ is the so-called {\it Bergman kernel}.
It is defined by the following absolutely and uniformly convergent series:
\begin{equation}
B_{\hbar}(z,\bar{v})=\sum_{k=1}^{\infty}f_k(z)\overline{f_k(v)},
\end{equation}
for an arbitrary orthonormal basis $\{f_k\}$ in $\diff{F}_{\hbar}$ with
respect to the scalar product~\eqref{sc-prod}.

Representation~\eqref{projector} shows that holomorphic functions
$\Phi_{\bar{v}}(z)=B_{\hbar}(z,\bar{v})$, parameterized by $v\in\CC^n$, form
a complete system in $\diff{F}_{\hbar}$. For every bounded operator $F$ in
$\diff{F}_{\hbar}$ define its {\it covariant symbol}~\cite{berezin2} with
respect to the complete system $\{\Phi_{\bar{v}}\}_{v\in\CC^n}$ as the
restriction to the diagonal $v=z$ of the following function
\begin{equation}
f(z,\bar{v})=\frac{(F\Phi_{\bar{v}},\Phi_{\bar{z}})}
{(\Phi_{\bar{v}},\Phi_{\bar{z}})},
\end{equation}
defined on $\CC^n\times\CC^n$. The function $f$ is holomorphic in $z$ and
anti-holomorphic in $v$.

Denote by $\diff{A}_{\hbar}$ the linear space of covariant symbols of bounded
operators in $\diff{F}_{\hbar}$; there is a natural inclusion
$\diff{A}_{\hbar}\subset\diff{A}$. Then, according to~\cite{berezin}, we have
the following result.

\begin{lemma} The space of covariant symbols $\diff{A}_{\hbar}$ is a
$\CC$-algebra (an associative algebra over $\CC$ with unit), with the $\star$ -
product given by the following integral
\begin{equation} \label{integral1}
(f_1\star_{\hbar} f_2)(z,\bar{z})=\int_{\CC^n}f_1(z,\bar{v})f_2(v,\bar{z})
\frac{B_{\hbar}(z,\bar{v})B_{\hbar}(v,\bar{z})}{B_{\hbar}(z,\bar{z})}
e^{-\Phi(v,\bar{v})/\hbar}d\mu_{\hbar}(v,\bar{v}),
\end{equation}
and with the unit given by the symbol of the identity operator --- function on
$\CC^n$, identically equal to $1$. Also, there exists a trace functional
$\tr:\diff{A}_{\hbar}^{0}\to\CC$, defined on symbols of bounded operators
of the trace class by the formula
\begin{equation}
\tr f=\int_{\CC^n}f(z,\bar{z})B_{\hbar}(z,\bar{z})e^{-\Phi(z,\bar{z})/\hbar}
d\mu_{\hbar}(z,\bar{z}),
\end{equation}
and it is cyclically invariant
\begin{equation*}
\tr(f_1\star_{\hbar} f_2)=\tr(f_2\star_{\hbar} f_1).
\end{equation*}
\end{lemma}

\begin{proof}
It is a straightforward computation; see~\cite{berezin, berezin2} for details.
The associativity of the $\star$ - product follows from the associativity of
the operator product. It can be also verified directly from the
integral~\eqref{integral1} using the Fubini theorem; all issues regarding
convergence, etc., are standard. The same applies to the cyclic trace property.
Finally, the property of the unit follows from the definition of the Bergman
kernel.
\end{proof}

\subsection{Contravariant symbols} A bounded operator $F$ in
$\diff{F}_{\hbar}$ is called a {\it Berezin-Toeplitz operator}, if $F=P\hat{f}$,
where $\hat{f}$ is a multiplication operator by a function
$\hat{f}\in\diff{H}_{\hbar}$, called the {\it contravariant symbol}
of $F$. Covariant symbol $f$ of the operator $F$ can be expressed through its
contravariant symbol $\hat{f}$ (if it exists) as follows~\cite{berezin2}
\begin{equation} \label{contra}
f(z,\bar{z})=\int_{\CC^n}\hat{f}(v,\bar{v})\frac{B_{\hbar}(z,\bar{v})
B_{\hbar}(v,\bar{z})}{B_{\hbar}(z,\bar{z})}e^{-\Phi(v,\bar{v})/\hbar}
d\mu_{\hbar} (v,\bar{v}).
\end{equation}
The proof is a straightforward computation. The correspondence $\hat{f}
\mapsto f$ defines a linear operator $I_{\hbar}: \diff{F}_{\hbar}\rightarrow
\diff{F}_{\hbar}$, $I_{\hbar}(\hat{f})=f$.
\begin{lemma} \label{Contra}If $F$ is a Berezin-Toeplitz
operator in $\diff{F}_{\hbar}$ with  contravariant symbol $\hat{f}$, then
for all $g\in\diff{F}_{\hbar}$ we have~\cite{berezin2}
\begin{equation} \label{contra2}
F g(z)=\int_{\CC^n}B_{\hbar}(z,\bar{v}) f(z,\bar{v})g(v)
e^{-\Phi(v,\bar{v})/\hbar}d\mu_{\hbar}(v,\bar{v}),
\end{equation}
where $f$ is the covariant symbol of $F$, given by the formula~\eqref{contra}.
Also, the product of Berezin-Toeplitz operators $F_1$ and $F_2$ with
contravariant symbols $\hat{f}_1$ and $\hat{f}_2$ can be expressed through
their covariant symbols $f_1$ and $f_2$ by the following formula:
\begin{equation} \label{contra3}
(F_1 F_2)g(z)=\int_{\CC^n}B_{\hbar}(z,\bar{v})(f_1\star
f_2)(z,\bar{v})e^{-\Phi(v,\bar{v})/\hbar}d\mu_{\hbar}(v,\bar{v}).
\end{equation}
\end{lemma}
\begin{proof} It is another straightforward computation, based on the
following expression of $\star$-product of covariant symbols through
corresponding contravariant ones,
\begin{gather} \label{contra4}
(f_1\star_{\hbar} f_2)(z,\bar{z})=\iint_{\CC^n\times\CC^n}\hat{f}_1(v,\bar{v})
\hat{f}_2(w,\bar{w})\frac{B_{\hbar}(z,\bar{v})B_{\hbar}(v,\bar{w})
B_{\hbar}(w,\bar{z})}{B_{\hbar}(z,\bar{z})} \\ \times
e^{-(\Phi(v,\bar{v})+\Phi(w,\bar{w}))/\hbar}d\mu_{\hbar}(v,\bar{v})
d\mu_{\hbar}(w,\bar{w}). \nonumber
\end{gather}
\end{proof}
If the linear operator $I_{\hbar}$ were one-to-one and onto,
equation~\eqref{contra4} would allow to introduce the product
$\hat{\star}_{\hbar}$ of contravariant symbols as follows
\begin{equation} \label{contra-star}
\hat{f}_1\hat{\star}_{\hbar}\hat{f}_2=I_{\hbar}^{-1}(I_{\hbar}(\hat{f}_1)
\star_{\hbar} I_{\hbar}(\hat{f}_2)).
\end{equation}
In general~\cite{berezin, berezin2}, however, the operator $I_{\hbar}$ is not
onto, since there exist bounded operators which do not have contravariant
symbols.

\subsection{The correspondence principle} To consider the family of
associative algebras $\{\diff{A}_{\hbar}\}_{\hbar>0}$ as a quantization of
the Poisson algebra of classical observables $\diff{A}=C^{\infty}(\CC^n)$,
the ``correspondence principle''
\begin{equation} \label{principle}
\lim_{\hbar\rightarrow 0}\frac{2\sqrt{-1}}{\hbar}(f_1\star_{\hbar}
f_2-f_2\star_{\hbar} f_1) =\{f_1,f_2\},
\end{equation}
should be valid for "suitable" continuous families $f_i(z,\bar{z}; \hbar)
\in\diff{A}_{\hbar}\subset\diff{A}$, where in the right-hand side
$f_i(z,\bar{z})=f_i(z,\bar{z};\hbar)\mid_{\hbar=0}\,\in\diff{A},~i=1,2$.

In order to verify~\eqref{principle} and to introduce these families of
symbols, Berezin~\cite{berezin} made several assumptions.

First, he assumed that K\"{a}hler potential $\Phi(z,\bar{z})$ admits an
analytic continuation $\Phi(z,\bar{v})$ to $\CC^n\times\CC^n$. Setting
\begin{equation}
B_{\hbar}(z,\bar{v})=e^{\Phi(z,\bar{v})/\hbar}e_{\hbar}(z,\bar{v}),
\end{equation}
one can rewrite the integral \eqref{integral1} as follows
\begin{equation} \label{integral2}
(f_1\star_{\hbar} f_2)(z,\bar{z})=\int_{\CC^n}f_1(z,\bar{v})f_2(v,\bar{z})
\frac{e_{\hbar}(z,\bar{v})e_{\hbar}(v,\bar{z})}{e_{\hbar}(z,\bar{z})}
e^{\phi(z, \bar{z};v, \bar{v})/\hbar}d\mu_{\hbar}(v,\bar{v}),
\end{equation}
where $\phi(z,\bar{z};v,\bar{v})=\Phi(z,\bar{v})+\Phi(v,\bar{z})-
\Phi(z,\bar{z}) -\Phi(v,\bar{v})$ is the so-called {\it Calabi diastatic
function} of the K\"{a}hler form $\omega$.

\begin{remark} \label{Calabi}
When K\"{a}hler form $\omega$ is real-analytic, the Calabi function is globally
defined in some neighborhood of the diagonal in $\CC^n\times\CC^n$. In general,
for an arbitrary K\"{a}hler manifold with a real-analytic K\"{a}hler form, the
Calabi function is defined in some neighborhood of the diagonal in $M\times M$.
\end{remark}

The exponential factor in integral~\eqref{integral2}  has a critical point at
$v=z$, so that the Laplace method is applicable (for smooth functions $f_1$ and
$f_2$), provided that the asymptotic behavior of $e_{\hbar}$ as $\hbar\to 0$
is also known. In~\cite{berezin}, Berezin did not consider this problem, which
we address and solve in section 4. Instead, Berezin made the  second assumption
(hypothesis A-D in~\cite{berezin}), that $e_{\hbar}(z,\bar{v})=1$ for all
$z,v\in\CC^n$.

This condition is very restrictive. It implies (see~\cite{berezin}, Theorem
2.5) that $\log\det H $ is harmonic, so that Berezin's quantization works only
for very special K\"{a}hler manifolds, essentially, for the flag varieties and their
generalizations. The corresponding construction depends on a particular choice of the
K\"{a}hler potential and is coordinate dependent.

\section{Non-normalized $\star$ - product for K\"{a}hler manifolds}

\subsection{Non-normalized $\star$ - product for $\CC^n$}  Let
$\Phi(z,\bar{v})$ be analytical continuation of $\Phi(z,\bar{z})$ to
$\CC^n\times
\CC^n$ and let $\phi(z,\bar{z};v,\bar{v})=\Phi(z,\bar{v})+
\Phi(v,\bar{z})-\Phi(z,\bar{z})-\Phi(v,\bar{v})$ be its Calabi function. The
following
statement holds.
\begin{lemma} \label{expansion}
The point $v=z$ is a critical point for the Calabi function $\phi(z,\bar{z};
v,\bar{v})$ considered as a function of $v~\text{and}~\bar{v}$. One has the
following expansion
\begin{equation*}
\phi(z,\bar{z}; v,\bar{v})=-\sum_{|I|+|J|>1}\frac{1}{I!J!}\Phi_{I\bar{J}}
(v-z)^I(\bar{v}-\bar{z})^J,
\end{equation*}
with the standard notation for the multi-indices $I=(i_1,\ldots,i_n),
J=(j_1,\ldots,j_n)$,
\begin{equation*}
\begin{split}
|I| & = i_1+\cdots+i_n,~|J|=j_1+\cdots+j_n,~I!=i_1!\cdots i_n!,~J!=j_1!\cdots
j_n!, \\
(v-z)^I & = (v^{1}-z^{1})^{i_1}\cdots (v^{n}-z^{n})^{i_n},~(\bar{v}-\bar{z})^J
=
(\bar{v}^{1}- \bar{z}^{1})^{j_1}\cdots(\bar{v}^{n}-\bar{z}^{n})^{j_n}.
\end{split}
\end{equation*}
Here
\begin{equation*}
\Phi_{I\bar{J}}=\del^I\delb^J\Phi=
\del^{I^{\prime}}\delb^{J^{\prime}}h_{i_1\bar{j_1}},
\end{equation*}
where $I^{\prime}=I\setminus\{i_1\},~J^{\prime}=J\setminus\{j_1\}$ and
$\del^I=\del_{1}^{i_1}\cdots \del_n^{i_n},~\delb^J=\delb_1^{j_1}\cdots\
\delb_n^{j_n}$.
In particular, $\Phi_{i\bar{j}}=\del_i\delb_j\Phi=h_{i\bar{j}}$.
\end{lemma}

This lemma shows that coefficients of the Taylor expansion of $\phi$ at the
critical point $v=z$ depend only on the K\"{a}hler metric. The following
definition is thus justified.

\begin{definition} The non-normalized $\star$ - product is a $\CC[[h]]$-linear
map
\begin{equation*}
\bullet:\diff{A}[[\hbar]]\otimes\diff{A}[[\hbar]]\rightarrow\diff{A}[[\hbar]],
\end{equation*}
defined by the following formal integral
\begin{equation} \label{nn}
(f_1\bullet f_2)(z,\bar{z})=\int_{\CC^n}f_1(z,\bar{v})f_2(v,\bar{z})
e^{\phi(z,\bar{z}; v,\bar{v})/\hbar}d\mu_{\hbar}(v,\bar{v}),
\end{equation}
i.e.~the formal power series in $\hbar$ that arises from the Laplace expansion
at the critical point $v=z$.
\end{definition}

\begin{remark} \label{formalint} In the formal integral~\eqref{nn} there is no
need to assume that functions $f_1$ and $f_2$ admit analytic continuation into
$\CC^n\times\CC^n$. The definition of the $\bullet$ - product uses the formal
expansion around the diagonal in $\CC^n\times\CC^n$ and depends only on the
values of $f_1$ and $f_2$ and of their partial derivatives on the diagonal.
Specifically, the integral can be written as the following formal power series
in $\epsilon^2=\hbar$, \begin{equation} \label{nn-1}
(f_1\bullet f_2)(z,\bar{z})=\int_{\CC^n}f_1(z,\bar{z}+\epsilon\bar{y})
f_2(z+\epsilon y, \bar{z})e^{\phi(z,\bar{z};z+\epsilon y,\bar{z}+
\epsilon\bar{y})/\epsilon^2} d\mu(y,\bar{y}).
\end{equation}
where $\phi(z,\bar{z};z+\epsilon y,\bar{z}+\epsilon\bar{y})$ is the formal
power series in $\epsilon$ determined by lemma~\ref{expansion}, $f_1(z,\bar{z}+
\epsilon\bar{y})$ and $f_2(z+\epsilon y,\bar{z})$, are given by their
formal power series expansions, $d\mu(y,\bar{y})=(\omega/\pi)^n$ and Gaussian
integrations have been performed (see~\eqref{gauss} in the next subsection).
Equivalently, the formal integral~\eqref{nn} can be defined by the asymptotic
of the distribution $\exp(\phi(z,\bar{z};v,\bar{v}) /\hbar)$ as $\hbar
\rightarrow 0$.
\end{remark}

\begin{lemma} \label{ass} The product $\bullet$ is associative,
\begin{equation*}
(f_1\bullet f_2)\bullet f_3=f_1\bullet (f_2\bullet f_3),
\end{equation*}
for all $f_1, f_2, f_3\in\diff{A}$.
\end{lemma}

\begin{proof}
It follows from the definition of the Calabi function that
\begin{equation*}
\phi(z,\bar{w},v,\bar{v})+\phi(z,\bar{z},w,\bar{w})=
\phi(v,\bar{z},w,\bar{w})+\phi(z,\bar{z},v,\bar{v}).
\end{equation*}
Therefore, both sides of the associativity equation are equal to the same
formal integral
\begin{equation*}
\iint_{\CC^n\times\CC^n}f_1(z,\bar{v})f_2(v,\bar{w})f_3(w,\bar{z})
e^{(\phi(z,\bar{z};v,\bar{v})+ \phi(v,\bar{z};w,\bar{w}))/\hbar}
d\mu_{\hbar}(v,\bar{v}) d\mu_{\hbar}(w,\bar{w}),
\end{equation*}
provided that one has the formal analog of the Fubini theorem: ``double formal
integral equals to the repeated formal integral''.

The latter statement can be proved as follows. First, it is clearly valid for
convergent integrals, provided that all sufficient conditions for the uniform
convergence are satisfied. These conditions are in the form of inequalities
for the coefficients $\Phi_{I\bar{J}}$ and for the Taylor coefficients of
$f_1, f_2, f_3$. Since the associativity equation can be written as a system
of polynomial equations for these coefficients, these equations hold for
arbitrary values of the coefficients involved --- ``the nonessential property
of algebraic identities" principle.

\end{proof}

\subsection{Feynman graphs and coefficients of non-normalized $\star$ -
product for $\CC^n$} Here we represent the formal power series for the
$\bullet$ - product using Feynman graphs (see   any textbook on
quantum field theory; for a mathematical exposition see, e.g.~\cite{kazhdan}).

Let $\G$ be a graph, i.e.~a finite one-dimensional simplicial complex with
oriented edges. The graph $\G$ is {\it locally oriented} if an
enumeration of identically oriented edges is fixed around each vertex.

Consider a locally oriented graph $\G$ equipped with the following data.
\begin{itemize}
\item[(i)] For each edge $e\in\G_1$ of $\G$ there are given a finite-dimensional
vector space $V_e$ over $\CC$, its dual space $V_e^{\prime}$, and a linear map
$L_e:V_e^{\prime}\rightarrow V_e$.
\item[(ii)] For each vertex $v\in\G_0$ of $\G$ there is given a linear map
\begin{equation*}
\g_v: V_{e_{1v}^{-}}\otimes\cdots\otimes V_{e_{nv}^{-}}\rightarrow
V_{e_{1v}^{+}}^{\prime}\otimes\cdots\otimes V_{e_{mv}^{+}}^{\prime},
\end{equation*}
equivariant with respect to the symmetric group, which acts by permuting the
factors in the tensor product. Here $e_{iv}^{-},~i=1,\dots, n$, are
incoming edges, and $e_{jv}^{+},~j=1, \dots, m$, are outgoing edges in the
star $\Gamma(v)$ of the vertex $v$. Also, it is assumed that the tensor product
of vector spaces parameterized by the empty set is $\CC$.
\end{itemize}
To define the partition function of a locally oriented graph $\G$
with data, choose for every edge $e\in\G_1$ a basis
$\{\e_{i_e}\}^{\dim V_e}_{i_e=1}$ in the vector space $V_e$ and let
$\{\e^{i_e}\}^{\dim V_e}_{i_e=1}$ be the dual basis in
$V_{e}^{\prime}$. Denote by $(\g_v)_{I_{-}I_{+}},~(L_e)^{k_e l_e}$, where
$I_{-}=\{i_{e^{-}_{1v}},\dots, i_{e^{-}_{nv}}\}$ and
$I_{+}=\{i_{e^{+}_{1v}},\dots,i_{e^{+}_{mv}}\}$, the corresponding matrix elements
of the linear maps $\g_v$ and $L_e$, and set
$N=h_1(\G)\stackrel{def}{=}\#(\G_1)$. Then the partition function of $\G$ is given
by the following expression,
\begin{equation*}
W_{\G}=W_{\G}(\{\g_v\},\{V_e\},\{L_e\})\stackrel{def}{=}
\sum_{i_{e_1}=1}^{\dim V_{e_1}}\cdots\sum_{i_{e_N}=1}^{\dim V_{e_N}}
\prod_{v\in\G_0} (\g_v)_{I_{-}I_{+}}\prod_{e\in\G_1} (L_e)^{k_e l_e},
\end{equation*}
where $k_e=i_{e^{+}_{pv}}$ for the unique $v\in\G_0$ such that $e=e^{+}_{pv}$
and $l_e=i_{e^{-}_{qv}}$ for the unique $v\in\G_0$ such that $e=e^{-}_{qv}$.

Due to the equivariance property of maps $\g_v$, the partition function $W_{\G}$
does not depend on the ordering of identically oriented edges around each
vertex of the graph $\G$. It also does not depend on the choice of bases in
vector spaces $V_e$, and is an invariant of the locally oriented graph $\G$
with data.

The partition functions of graphs with data are building blocks for expressing the
asymptotics of integrals through Feynman diagrams. Specifically, consider the
formal integral~\eqref{nn} which, by remark~\ref{formalint},
we rewrite as follows,
\begin{equation}
\pi^{-n} \int_{\CC^n}f_1(z,\bar z +\e\bar y) f_2(z+\e y, \bar z)e^{-(Hy,y)+
V(z,\bar z;y,\bar y;\hbar)} \prod^n_{i=1} \frac{|dy^i\wedge d\bar y^i|}{2}.
\end{equation}
Here $\e=\sqrt{\hbar}$ and
\begin{equation*}
(Hy,y)=\sum^n_{i,j=1}h_{i\bar j}(z,\bar z)y^i\bar{y}^j,
\end{equation*}
where $H=(h_{i\bar j})$ is the matrix of the K\"{a}hler metric, and
\begin{eqnarray*}
V(z,\bar z;y,\bar y;\hbar)&=&-\sum_{|I|+|J|>2}
\frac{\hbar^{(|I|+|J|-2)/2}}{I!J!}\Phi_{I\bar{J}} (z,\bar z)y^I\bar{y}^J \\
& &\mbox{}+\sum_{|I|,|J|\geq 0}\frac{\hbar^{(|I|+|J|)/2}}{I!J!}\Psi_{I\bar{J}}
(z, \bar z)y^I \bar{y}^J,
\end{eqnarray*}
where we are using notations of Lemma 3.1 and
\begin{equation*}
\Psi_{I\bar{J}}(z,\bar z)=\del^I\delb^J\log\det H(z,\bar z).
\end{equation*}

The evaluation of the formal integral~\eqref{nn} is based on the fundamental
formula of Gaussian integration
\begin{equation} \label{gauss}
\int_{\CC^n}y^I\bar{y}^J e^{-(Hy,y)}\prod_{i=1}^n\frac{|dy^i\wedge d
\bar{y}^i|}{2}=\frac{\pi^n}{\det H}(\del^I\delb^J e^{(H^{-1}z,z)})\mid_{z=
\bar z=0}.
\end{equation}
Namely, expand $f_1(z,\bar z+\e\bar y)f_2(z+\e y,\bar z)\exp V(z,\bar z;y,
\bar y;\hbar)$ in power series in $\e$ and evaluate each term
using~\eqref{gauss}. As a result, one gets an expression of~\eqref{nn} as a
formal power series in $\hbar=\e^2$ with coefficients given as sums over a subset
$\diff{G}$ of the set of all locally oriented graphs.

This subset $\diff{G}$ is characterized by the two conditions. First, for every
$\G\in\diff{G}$, the set $\G_0$ always contains two special elements: a
vertex $R$ with only incoming edges, and a vertex $L$ with only outgoing
edges (see Fig. 1).

\begin{figure}[!h]
\centering
\epsfig{figure=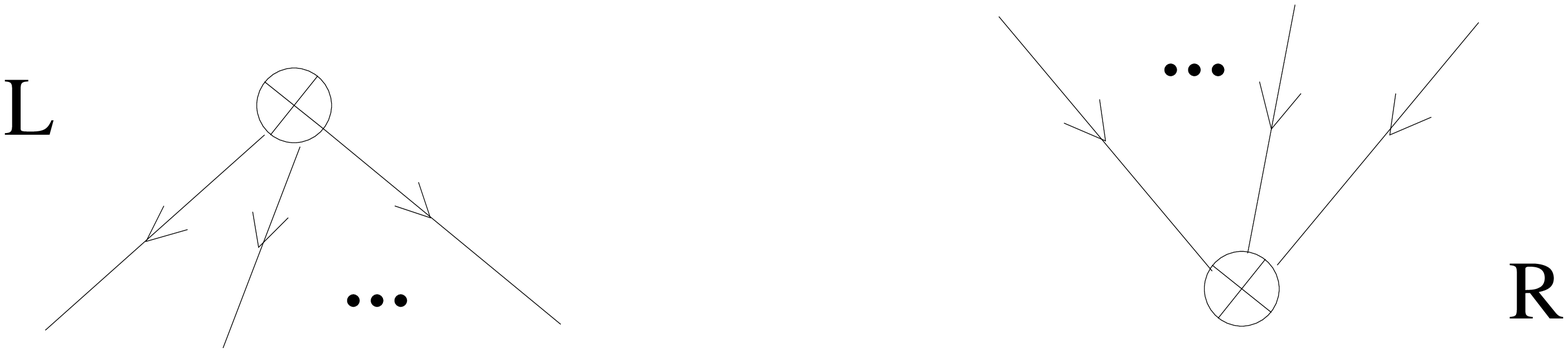,height=2cm,angle=0}
\caption{}
\end{figure}

Second, the complement of the set $\{L,R\}$ in $\G_0$ consists of vertices of
two types, ``solid" and ``hollow", either one or both of these sets may be
empty (see Fig. 2).

\begin{figure}[!h]
\centering
\epsfig{figure=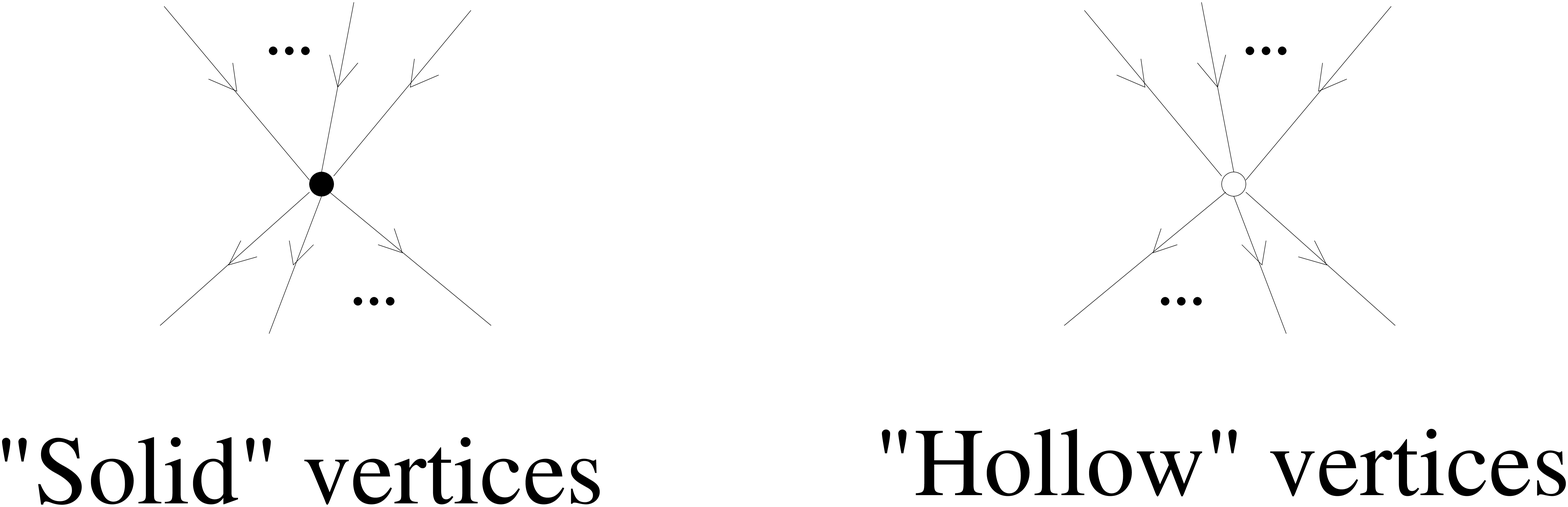,height=2cm,angle=0}
\caption{}
\end{figure}

For every $\G\in\diff{G}$ let $\Aut(\G)$ be the subgroup of the group of
automorphisms of $\G$ that maps vertices to the vertices of the same type and
preserves $L$ and $R$ vertices and let $|\Aut(\G)|=\#\Aut(\G)$. Also set
$\chi(\G)=\#(\G_1)-\#(\G_{0}^{s})$, where $\G_{0}^s\subset\G_0$ is the subset
of ``solid'' vertices.

The next theorem summarizes the description of the formal integral~\eqref{nn} in
terms of Feynman diagrams.

\begin{theorem} \label{fgraphs} The non-normalized $\star$ - product
$\bullet$ is given by
\begin{equation*}
f_1\bullet f_2 =\sum_{\G\in\diff{G}}
\frac{\hbar^{\chi(\G)}}{|\Aut(\G)|}D_{\G}(f_1,f_2),
\end{equation*}
where $D_{\G}(f_1,f_2)=W_{\Gamma}$ is the partition function of the graph
$\G\in\diff{G}$ with the following data.
\begin{itemize}
\item[1.] For each edge $e\in\G_1$ there are given a vector space
$V=\CC^n$ with the standard basis $\{\e_i\}$, its dual space $V^{\prime}=
\CC^n$ with the dual basis $\{\e^{\bar i}\}$ and a linear map $L:
V^\prime\rightarrow V$, defined by $K(\e^{\bar j})=\sum_{i=1}^n h^{i\bar{j}}
(z,\bar z)\e_i$, where $h^{i\bar{j}}$ are elements of the matrix $H^{-1}$.
\item[2.] For each ``solid" vertex $v_s$, with $n$ incoming
and $m$ outgoing edges, there is given a linear map $\g_{v_s}: V^{\ot n}
\rightarrow (V^\prime)^{\ot m}$, defined by the matrix with elements
$-\Phi_{I\bar{J}}(z,\bar z)$.
\item[3.] For each ``hollow" vertex $v_h$, with $n$ incoming
and $m$ outgoing edges, there is given a linear map $\g_{v_h}: V^{\ot n}
\rightarrow (V^{\prime})^{\ot m}$, defined by the matrix with elements
$ \Psi_{I\bar{J}}(z,\bar z)$.
\item[4.] For the special vertex $R$, with $n$ incoming edges, there is given a
linear map $\g_R: V^{\ot n}\rightarrow \CC$, defined by the matrix
with elements $f_{2 I}(z\bar z)=\del^I f_2(z,\bar z)$.
\item[5.] For the special vertex $L$, with $m$ outgoing
edges, there is given a linear map $\g_L: \CC\rightarrow (V^{\prime})^{\ot n}$
defined by the matrix with elements $f_{1 \bar J}(z,\bar z)=
\delb^J f_1(z,\bar z)$.
\end{itemize}
\end{theorem}

\begin{remark} \label{relations}
Elements of linear maps assigned to solid and hollow vertices are both
expressed through the K\"{a}hler metric and, therefore, are related. These
relations are obtained from the fundamental one
\begin{equation*}
\p_i\log\det H=\tr (H^{-1}\p_i H),
\end{equation*}
which is represented graphically in Figure 3, by successive differentiation.

\begin{figure}[!h]
\centering
\epsfig{figure=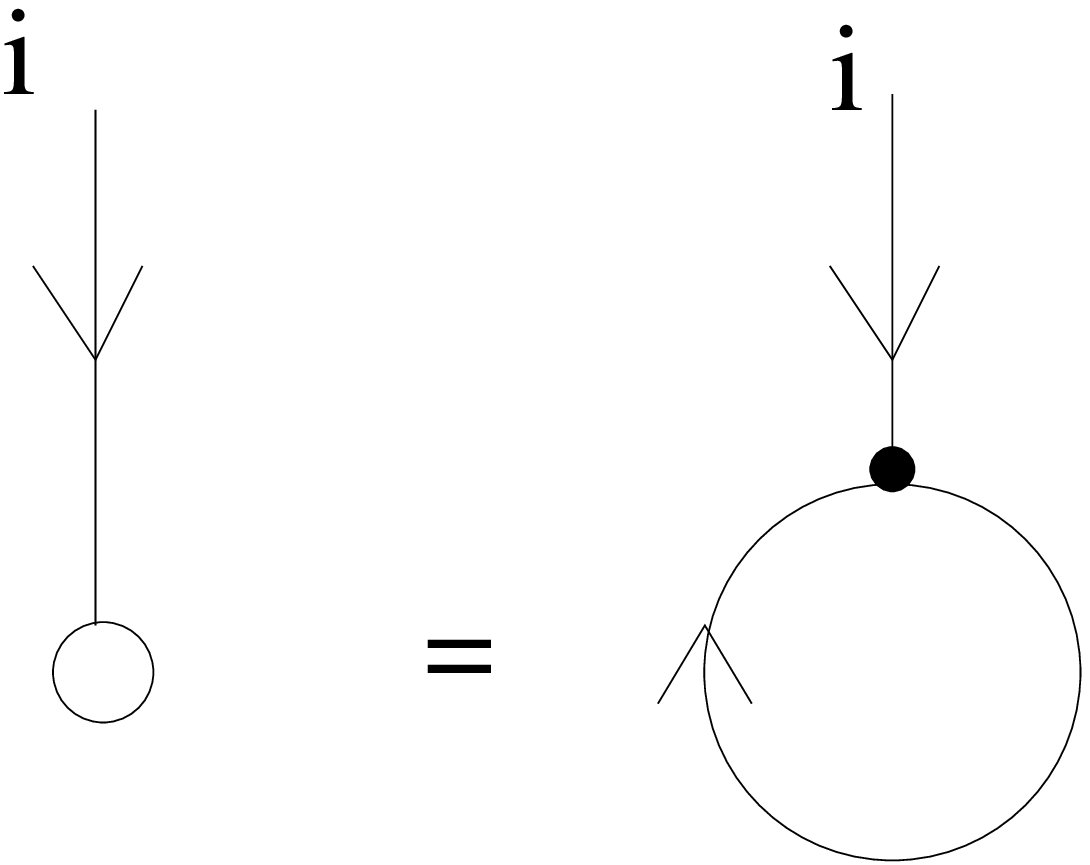,height=2cm,angle=0}
\caption{}
\end{figure}

Typical examples are
\begin{equation*}
\p^{2}_{ij}\log\det H=\tr(H^{-1}\p^{2}_{ij}H)- \tr(H^{-1}\p_i H H^{-1}\p_j H),
\end{equation*}
depicted as equality of weights of graphs in Figure 4, and
\begin{eqnarray*}
\p^{3}_{ijk}\log\det H & = & \tr(H^{-1}\p_i H H^{-1}\p_j H H^{-1}\p_k H)\\
&-& \tr(H^{-1}\p^2_{ij} H H^{-1} \p_k H ) +\tr(H^{-1}\p_j HH^{-1}\p_k HH^{-1}
\p_i H)\\
&-&\tr(H^{-1}\p^{2}_{ik} HH^{-1}\p_j H)-\tr(H^{-1}\p_i HH^{-1}\p^{2}_{jk} H) \\
&+& \tr(H^{-1}\p^{3}_{ijk} H),
\end{eqnarray*}
presented in Figure 5.
\end{remark}

\begin{figure}[!h]
\centering
\epsfig{figure=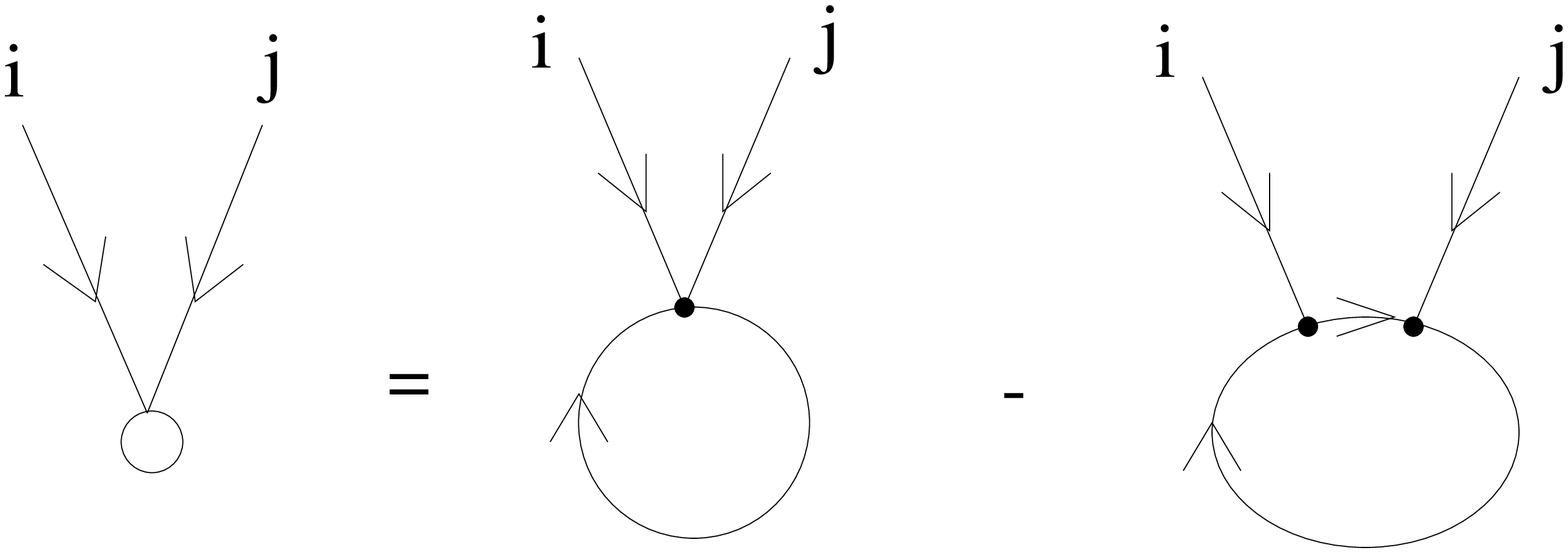,height=2cm,angle=0}
\caption{}
\end{figure}

\begin{figure}[!h]
\centering
\epsfig{figure=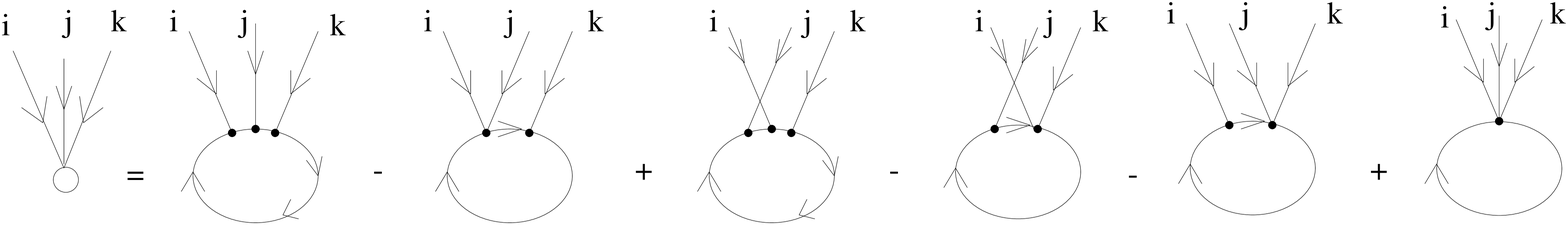,height=2cm,angle=0}
\caption{}
\end{figure}

\begin{lemma} The first two terms of the product
$f_1\bullet f_2$ are given by
\begin{eqnarray}
 f_1\bullet f_2&=&f_1f_2+\hbar(Af_1f_2+\sum_{i,j=1}^n h^{i\bar j}\delb_j f_1
\del_i f_2) \nonumber\\
&+&\frac{\hbar^2}{2}(\sum_{i,j=1}^n\sum_{k,l=1}^n
(h^{i\bar{j}} h^{k\bar{l}} \delb^{2}_{jl}f_1 \del^{2}_{ik}
f_2\nonumber\\
&+& h^{i\bar{j}}\del_{i}h^{k\bar{l}}\delb^{2}_{jl}f_1 \del_{k} f_2
+ h^{i\bar{j}}\delb_{j}h^{k\bar{l}}\delb_{l}f_1 \del^{2}_{ik} f_2)
\nonumber \\
&+&\mbox{} \sum_{i,j=1}^n h^{i\bar j}(\delb_j(Af_1)\del_i f_2 + \delb_j
f_1\del_i(Af_2) \\
&+&\mbox{} \sum_{i,j=1}^n Ah^{i\bar j}\delb_j f_1\del_i f_2
+2Df_1f_2)+O(\hbar^3).\nonumber
\end{eqnarray}
Here
\begin{equation*}
A=\frac{1}{2}\sum_{i,j=1}^n h^{i\bar j}\del_i\delb_j\log\det H,
\end{equation*}
and the term $\hbar^2 D f_1f_2$ is the contribution from the vacuum diagrams, i.e.,
graphs with no edges adjacent to special vertices $L$ and $R$.
\end{lemma}

\begin{proof} It is an easy computation, based on theorem~\ref{fgraphs}
and remark~\ref{relations}.
\end{proof}
Diagrammatically, this expression is given in Figure 6.

\begin{figure}[!h]
\centering
\epsfig{figure=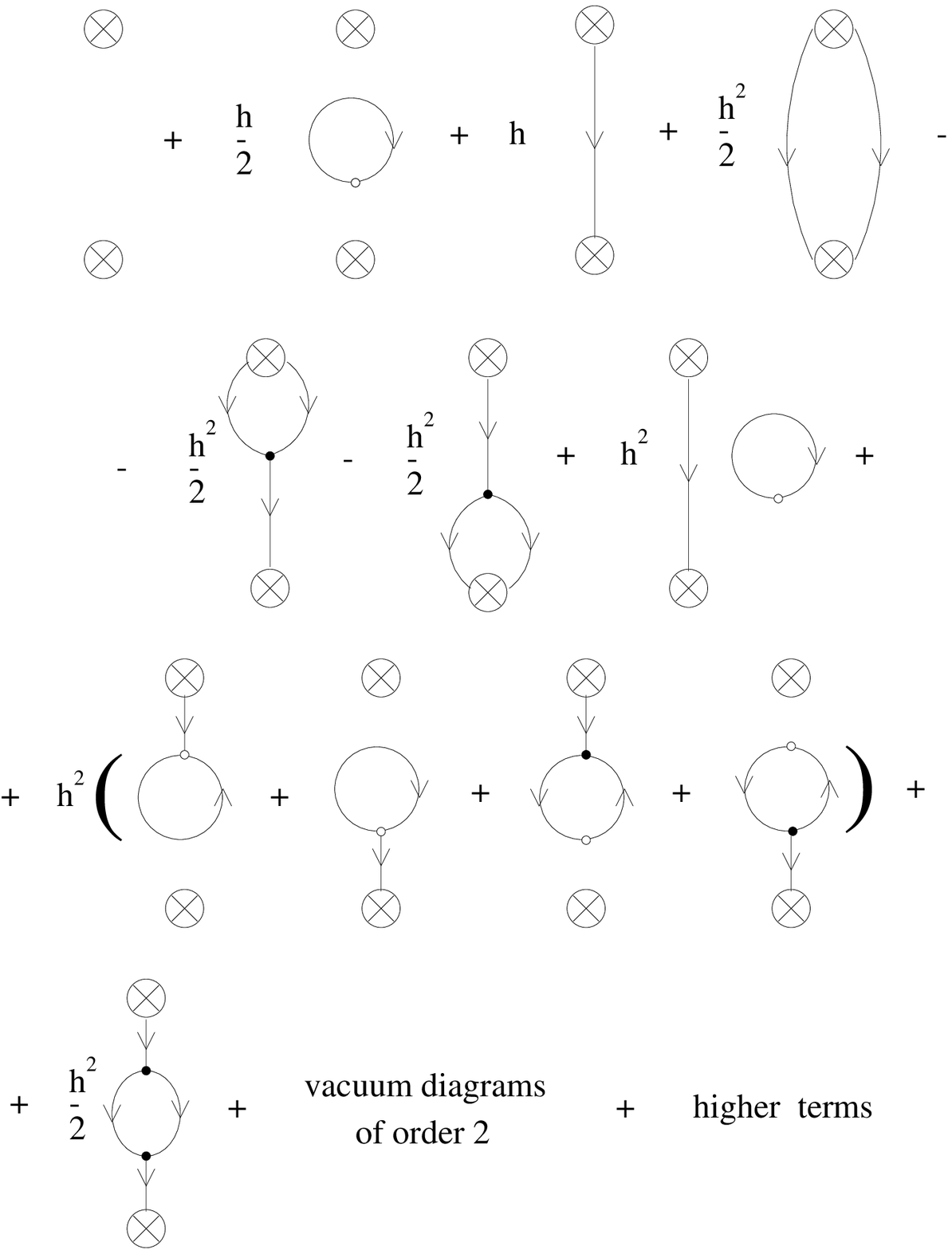,height=14cm,angle=0}
\caption{}
\end{figure}

\begin{corollary}
The product $\bullet$ gives a deformation quantization of the symplectic
manifold $(\CC^n, \omega)$ with the property $\overline{f_1\bullet
f_2}=\bar{f_2}\bullet\bar{f_1}$, where $\bar{f}(z,\bar{z})=
\overline{f(\bar{z},z)}$.
\end{corollary}

\begin{proof}
As it follows from the previous lemma,
\begin{equation}
\begin{split}
f_1\bullet f_2 &= f_1f_2 + O(\hbar), \\ f_1\bullet f_2-f_2\bullet f_1 &=
\frac{\hbar}{2\sqrt{-1}}\{f_1,f_2\}+O(\hbar^2),
\end{split}
\end{equation}
so that the correspondence principle is valid. The symmetry with respect to
the complex conjugation follows from the definition of the $\bullet$ -
product.
\end{proof}

\begin{remark} \label{fermions}
Another way to represent the formal integral~\eqref{nn} in terms
of Feynman diagrams is to rewrite it as a formal integral over the
supermanifold $\Pi T\CC^n$ --- the total space of the tangent bundle over
$\CC^n$ with reversed parity of the fibers. Namely, introducing Gaussian
integration over the Grassmann variables $\xi^i, \bar\xi^i,~i=1, \cdots,n$,
we get
\begin{eqnarray*}
& &\mbox{} \int_{\CC^n} f_1(z,\bar v)f_2(v, \bar z)e^{\phi(z,\bar z;v,\bar v)/
\hbar} \omega^n \\ &=&\int_{\Pi T\CC^n} f_1(v,\bar z)f_2(v, \bar z)
e^{\{\phi(z, \bar z; v,\bar v)/\hbar +(H(v, \bar v)\xi,\xi)\}}\prod_{i=1}^{n}
\frac{|dz^i\wedge d\bar z^i|}{2}d\xi^i d\bar\xi^i,
\end{eqnarray*}
where $(H\xi,\xi)=\sum_{i,j=1}^n h_{i \bar j}\xi^i\bar \xi^j$. Expanding into
power series in $\e=\sqrt{\hbar}$ through $v=z+\e y,~\bar v=\bar z+\e\bar y$,
and using the fundamental formula of Gaussian integration, we obtain a
representation of $f_1 \bullet f_2$ in terms of Feynman diagrams with
"bosonic" and "fermionic" vertices.
\end{remark}

\subsection{The trace}
On the subspace $\diff{A}_c\subset\diff{A}$ of compactly supported functions
consider the following $\CC[[\hbar]]$ - linear functional
$\tau:\diff{A}_c[[\hbar]]\rightarrow\CC[\hbar^{-1},\hbar]]$
\begin{equation}
\tau(f)=\int_{\CC^n}f(z,\bar{z})d\mu_{\hbar}(z,\bar{z}).
\end{equation}

\begin{proposition}
The functional $\tau$ is cyclically invariant:
\begin{equation*}
\tau(f_1\bullet f_2)=\tau(f_2\bullet f_1).
\end{equation*}
\end{proposition}

\begin{proof} It follows from the definitions that $\tau(f\bullet g)
\in\CC[\hbar^{-1},\hbar]]$ is given by the Laplace expansion around $v=z$ of
the following formal double integral:
\begin{equation*}
\iint_{\CC^n\times\CC^n}f_1(z,\bar{v})f_2(v,\bar{z})e^{\phi(z,\bar{z};v,
\bar{v})/\hbar} d\mu_{\hbar}(v,\bar{v})d\mu_{\hbar}(z,\bar{z}),
\end{equation*}
which is symmetric with respect to $v$ and $z$.
\end{proof}

\begin{remark} In terminology of~\cite{cfs}, the existence of the cyclic
trace functional implies that the $\bullet$ is a strongly closed $\star$ -
product.
\end{remark}

\begin{remark} Instead of the subspace $\diff{A}_c$ of compactly supported
functions on $\CC^n$, one could consider a larger subspace for which the
functional $\tau$ makes sense. However, this subspace would depend on the
global properties of the K\"{a}hler potential $\Phi$. The subspace $\diff{A}_c$
is the maximal in $\diff{A}$ with the property that $\tau$ is well-defined
for an arbitrary K\"{a}hler potential.
\end{remark}

\subsection{The non-normalized $*$-product for K\"{a}hler manifolds}

Potentials for the K\"{a}hler metric $ds^2$ on $\CC^n$ form a family
$\Phi(z,\bar{z})+f(z)+\overline{f(z)}$, parameterized by entire functions $f$.
However, the Calabi function $\phi(z,\bar{z};v,\bar{v})$   does not depend on the
choice of $f$ and is a well-defined function in a neighborhood of the
diagonal in $\CC^n\times\CC^n$ (cf.~remark~\ref{Calabi}).

The following theorem states that the $\bullet$ - product is ``functorial'' with
respect to the holomorphic change of coordinates.

\begin{theorem}
Let $w=g(z)$ be a holomorphic change of coordinates in a neighborhood
$U\subset\CC^n$ of a point $z_0\in\CC^n$. Then
\begin{equation}
(f_1\bullet f_2)\circ g|_U=((f_1\circ g)\bullet (f_2\circ g))|_U,
\end{equation}
where $(f\circ g)(z,\bar{z})=f(g(z),\overline{g(z)})$.
\end{theorem}

\begin{proof}
The integral~\eqref{nn} is invariant with respect to a local change of
coordinates $z\mapsto w=g(z)$ since $\phi(z,\bar{z};v,\bar{v})$ is a
well-defined function on $U\times U$. The same applies to the power series
expansion of~\eqref{nn}, which is ``the only thing that actually exists''.
\end{proof}

\begin{corollary}
Formula~\eqref{nn} provides a deformation quantization of an arbitrary complex
symplectic manifold --- a complex manifold $M$ with a symplectic form
$\omega\in\Omega^{(1,1)}(M)$.
\end{corollary}

\begin{proof}
It is an immediate consequence of  theorem~\ref{fgraphs}. Consider an
open covering of $M$,
\begin{equation*}
M=\cup_{\alpha\in A} U_{\alpha},
\end{equation*}
with the property that for every $\alpha\in A$ the forms
$\omega_{\alpha}=\omega|_{U_{\alpha}}$ admit a representation~\eqref{Phi}-
\eqref{omega} with some K\"{a}hler potential $\Phi_{\alpha}$. On every $U_{\alpha}$
define the $\star$ - product $f_1\bullet f_2$ of $f_1,f_2\in
C^{\infty}(U_{\alpha})$ as before. Though K\"{a}hler potentials $\Phi_{\alpha}$ on
$U_{\alpha}$ do not necessarily agree on the intersections $U_{\alpha}\cap U_{\beta}$,
the corresponding Calabi functions $\phi_{\alpha}$ on $U_{\alpha}$, which exist for
real-analytic $\omega$, combine to a well-defined function in some neighborhood of the
diagonal in $M\times M$. Applying the previous theorem in the real-analytic case,
we conclude that $f_1\bullet f_2 \in C^{\infty}(M)[[\hbar]]$. Since the
coefficients of the $\star$ - product are polynomials in elements of $H$ and
$H^{-1}$ and their partial derivatives, the same conclusion holds for the smooth case.
\end{proof}

\section{Normalized $\star$ - product for K\"{a}hler manifolds}
\subsection{Normalized $\star$ - product for $\CC^n$}
Although the product $\bullet$ solves the problem of deformation quantization
for complex symplectic manifolds, it has two drawbacks. First, its power series
expansion is given by bidifferential operators whose coefficients contain
elements of the matrix $H^{-1}$ as well as elements of $H$. This implies that
the product $\bullet$ is not suitable for the deformation quantization of
complex Poisson manifolds, where the coefficients should depend on the
elements of the inverse matrix $H^{-1}$ only. Second, the product $\bullet$ does
not preserve the unit element $1$ in the commutative algebra $\diff{A}$.

The second problem can be easily resolved. It is well-known~\cite{ger2} in the
deformation theory of the associative unital algebras that the unit element for
the deformed product always exists. It easily follows from the fact that
Hochschild cohomology classes $H^{*}(A,A)$ for a commutative $\CC$-algebra $A$
can be represented by  normalized cocycles. Indeed, normalizing the
coefficients of the $\star$ - product $\bullet$ in $A[[\hbar]]$, one gets the
unit element for $\bullet$ as a formal power series expressed through these
coefficients. Thus in our case, the unit element is given by
\begin{equation}
e_{\hbar}(z,\bar{z})=1+\sum_{l=1}^{\infty}\hbar^l
e^{(l)}(z,\bar{z})\in\diff{A}[[\hbar]],
\end{equation}
and is uniquely characterized by the equations
\begin{equation}
e_{\hbar}\bullet f=f\bullet e_{\hbar}=f,
\end{equation}
for all $f\in\diff{A}$. Because of the symmetry property with respect to the
complex conjugation, only one of these equations is independent. Specifically,
for the coefficients $e_l$ one has the following system of
equations
\begin{equation}\label{unit}
\sum_{l=0}^{\infty}\hbar^l\int_{\CC^n}e^{(l)}(z,\bar{v})f(v,\bar{z})
e^{\phi(z,\bar{z}; v,\bar{v})/\hbar} d\mu_{\hbar}(v,\bar{v})=f(z,\bar{z}),
\end{equation}
for all $f\in\diff{A}$.

\begin{remark}
It is not difficult to find the first few coefficients of the unit
$e_{\hbar}$. Indeed, using lemma 3.7 and the notations there, we get
\begin{equation} \label{unit2}
e_{\hbar}=1-\hbar A+\hbar^2(A^2-D)+O(\hbar^3),
\end{equation}
so that $e_{\hbar}(z,\bar{z})$ is an invertible element in $\diff{A}[[\hbar]]$.
\end{remark}

Now we can introduce the normalized $\star$ - product by twisting the
non-normalized product $\bullet$ by the unit element $e_{\hbar}$.

\begin{definition} The normalized product $\star$ is given by
\begin{equation}
(f_1\star f_2)(z,\bar{z})=e^{-1}_{\hbar}(z,\bar{z})((f_1e_{\hbar})\bullet
(f_2e_{\hbar}))(z,\bar{z}),
\end{equation}
where $f_1e_{\hbar}$ and $f_2e_{\hbar}$ stand for the point-wise product in
$\diff{A}[[\hbar]]$ and $e^{-1}_{\hbar}\in\diff{A}[[\hbar]]$.
\end{definition}

Explicitly, the $\star$-product is given by the following formal integral
\begin{equation} \label{n}
(f_1\star f_2)(z,\bar{z})=\int_{\CC^n}f_1(z,\bar{v})f_2(v,\bar{z})
\frac{e_{\hbar}(z,\bar{v})e_{\hbar}(v,\bar{z})}{e_{\hbar}(z,\bar{z})}
e^{\phi(z,\bar{z};v,\bar{v})/\hbar}d\mu_{\hbar}(v,\bar{v}).
\end{equation}

\begin{proposition} The $\star$ - product \eqref{n} endows
$\diff{A}[[\hbar]]$ with a structure of a $\CC$-algebra, and has the property
\begin{equation}\label{strong} f_1\star
f_2=f_1f_2+\hbar\sum_{i,j=1}^nh^{i\bar{j}}\delb_j f_1 \del_i f_2+O(\hbar^2).
\end{equation}
Moreover, the linear functional $\tr: \diff{A}_c[[\hbar]] \rightarrow
\CC[\hbar^{-1},\hbar]]$
\begin{equation*}
\tr(f)=\int_{\CC^n}e_{\hbar}(z,\bar{z})f(z,\bar{z})d\mu_{\hbar}(z,\bar{z}),
\end{equation*}
is cyclically invariant:
\begin{equation*}
\tr(f_1\star f_2)=\tr(f_2\star f_1).
\end{equation*}
\end{proposition}

\begin{proof}
The property $f\star 1=1\star f=f$, for all $f\in\diff{A}$, follows from the
definition. The strong form~\eqref{strong} of the correspondence principle
follows from the computation of $e^{(1)}$, given in Remark 4.1. The proof
of the trace property is the same as in proposition 3.10.
\end{proof}

\begin{remark}
Using the previous remark, it is easy to find the first two terms of the
normalized $\star$ - product. They are expressed only through elements of the
inverse matrix $H^{-1}=(h^{i\bar j})$ --- ``polarized Poisson tensor'', and
are given by the following expression
\begin{eqnarray*}
f_1 \star f_2 &=& f_1 f_2 + \hbar\sum_{i,j=1}^n h^{i\bar j} \delb_j f_1
\del_i f_2 +\frac{\hbar^2}{2}\sum_{i,j=1}^n\sum_{k,l=1}^n (h^{i \bar j}h^{k
\bar l} \delb^{2}_{jl} f_1 \del^{2}_{ik} f_2 \\
& + & h^{i \bar j} \delb_l (h^{k \bar l}) \delb_j f_1 \del^{2}_{ik} f_2 +
\del_i(h^{i \bar j}) h^{k \bar l} \delb^{2}_{jl} f_1\del_k f_2 \\
&+& \del_i (h^{i \bar j})\delb_l(h^{k \bar l}) \delb_j f_1 \del_k f_2) +
O(\hbar^3).
\end{eqnarray*}
\end{remark}
\begin{remark} Though the first two terms in $\star$ - product~\eqref{n}
are expressed through the Poisson tensor only, this is no longer true (without
further restrictions on the K\"{a}hler metric) for the higher order terms.
Therefore, the $\star$ - product~\eqref{n} is not suitable for the deformation
quantization of a complex Poisson structure on $\CC^n$ --- the problem
solved by Kontsevich for arbitrary Poisson structures on $\RR^n$~\cite{kon}.
Still, we believe that our approach can be successfully applied to the
problem of deformation quantization of arbitrary complex Poisson structure on
$\CC^n$, if it is modified in the spirit of the supersymmetric approach to
the Duistermaat-Heckman localization formula (see, i.e.,~\cite{s-z}).
Specifically, consider representation~\eqref{nn} for the non-normalized
$\bullet$ - product in the remark~\ref{fermions}, and notice that the
"action" $S=\phi(z,\bar z;v, \bar v)/\hbar +(H\xi,\xi)$ is not invariant
under the supersymmetry transformation. Introduction of the unit into the
normalized $\star$ - product~\eqref{n} modifies the action in a way that
cancels ``vacuum diagrams'', but still it does not restore the supersymmetry.
We conjecture that the supersymmetric action can be obtained by an
appropriate deformation ${\rm mod} \hbar$ of the K\"{a}hler metric. Then,
according to the paradigm of the supersymmetric quantum field theory, the
corresponding Feynman diagrams will contain no loops and will depend only on
the Poisson tensor and its derivatives.
\end{remark}

\subsection{Formal Bergman kernel} Here we will show that the unit element
$e_{\hbar}(z,\bar{z})$ can be interpreted as the kernel of a projection operator.
Specifically, for any $f\in\diff{A}$ define the following formal integral
\begin{equation} \label{oint}
P_{\hbar}(f)(z)=\oint_{\CC^n}e_{\hbar}(z,\bar{v})f(v,\bar{v})
e^{(\Phi(z,\bar{v})-\Phi(v,\bar{v}))/\hbar}d\mu_{\hbar}(v,\bar{v})
\end{equation}
as the formal asymptotic expansion around the point $v=z$ obtained by repeated
"integration by parts". Namely, applying Stokes' formula $n$ times, reduce
the integration of $(n,n)$-form
\begin{equation*}
e^{(\Phi(z,\bar{v})-\Phi(v,\bar{v})/\hbar}dv^1\wedge\ldots\wedge dv^n\wedge
d\bar{v}^1\wedge\ldots\wedge d\bar{v}^n
\end{equation*}
over the $2n$-chain---direct product of $n$ annuli centered at $z^i$ to the
integral of $(n,0)$-form over the $n$-cycle---direct product of $n$ circles of
radii $\epsilon$ around $z^i,~i=1,\ldots,n$, and iterate this procedure.
The term-by-term evaluation of the limit $\epsilon\rightarrow 0$
results in a power series expansion in $\hbar$ (see Appendix for more
details). Equivalently, the formal integral $\oint$
can be defined as the asymptotic of the distribution $\exp(\Phi(z,\bar{v})-
\Phi(v,\bar{v}))/\hbar)$ as $\hbar\rightarrow 0$.

This procedure is well-defined, since
\begin{equation*}
\Phi(z,\bar{v})-\Phi(v,\bar{v})=\sum_{i=1}^n\del_i\Phi(z,\bar{v})(z^i-v^i)
+\text{higher order terms},
\end{equation*}
and the matrix $\del_i\delb_j\Phi(z,\bar{z})=h_{i\bar{j}}(z,\bar{z})$ is
non-degenerate. Coefficients of this power series expansion in $\hbar$ are given by
infinite sums, convergent only when function $f$ and K\"{a}hler potential
$\Phi$ are real-analytic. In the $C^{\infty}$-category, these coefficients should be
interpreted as infinite jets. By $\CC[[\hbar]]$-linearity, $P_{\hbar}$ extends to the
linear mapping $P_{\hbar}: \diff{A}[[\hbar]]\rightarrow \diff{J}[[\hbar]]$, where
$\diff{J}$, where $\diff{J}$ is the space of global sections of the bundle of formal
(with respect to $\hbar$) infinite jets over $\CC^n$.

\begin{proposition} The mapping $P_{\hbar}$ takes values $\diff{J}_{hol}$ --- the space
of global sections of of the bundle of formal holomorphic jets over $\CC^n$.
\end{proposition}
\begin{proof} It is clear that in the real-analytic category the integral (\ref{oint}),
when it converges, is a holomorphic function. For the asymptotics of the integral, this is
an algebraic property of coefficients. Therefore, the same is true for the formal power
series with the coefficients in infinite jets in the $C^{\infty}$-category.
\end{proof}

\begin{proposition} Operator $P_{\hbar}$ is a projector,
$P_{\hbar}^2=P_{\hbar}$.
\end{proposition}
\begin{proof} This follows from the identity between two types of formal
integrals
\begin{equation} \label{identity}
\int_{\CC^n}f(v,\bar{v})e^{\phi(z,\bar{z};v,\bar{v})/\hbar}d\mu_{\hbar}
(v,\bar{v})=\oint_{\CC^n}\tilde{f}_{z}(v,\bar{v})e^{(\Phi(z,\bar{v})-
\Phi(v,\bar{v}))/\hbar}d\mu_{\hbar}(v,\bar{v}),
\end{equation}
where $\tilde{f}_{z}(v,\bar{v})=f(v,\bar{v})e^{(\Phi(v,\bar{z})-
\Phi(z,\bar{z}))/ \hbar}$, proved in Appendix. Namely, for $f\in\diff{F}$ and
$z\in\CC^n$ set $g_z(v)=f(v)e^{-(\Phi(v,\bar{z})-\Phi(z,\bar{z})/\hbar}$.
Using the property of the unit and that $g_z(v)$ is holomorphic in $v$, we get
\begin{eqnarray*}
f(z)& = & g_z(z)=(e_{\hbar}\bullet g_z)(z,\bar{z}) \\
    & = & \int_{\CC^n}e_{\hbar}(z,\bar{v})g_z(v)e^{\phi(z,\bar{z};v,\bar{v}))/
\hbar} d\mu_{\hbar}(v,\bar{v}) \\
    & = & \oint_{\CC^n}
e_{\hbar}(z,\bar{v})f(v)e^{(\Phi(z,\bar{v})-\Phi(v,\bar{v}))/\hbar}
d\mu_{\hbar}(v,\bar{v}) \\
    & = & P_{\hbar}(f)(z).
\end{eqnarray*}
\end{proof}
\subsection{Contravariant quantization} Here we define the formal analog of
the mapping $I_{\hbar}$ from contravariant to covariant symbols as the
following formal integral
\begin{equation} \label{I}
I_{\hbar}(\hat{f})(z,\bar{z})=\int_{\CC^n}\hat{f}(v,\bar{v})\frac{e_{\hbar}
(z,\bar{v})e_{\hbar}(v,\bar{z})}{e_{\hbar}(z,\bar{z})}e^{\phi(z,\bar{z};
v,\bar{v})/\hbar}d\mu_{\hbar}(v,\bar{v}).
\end{equation}
By $\CC[[\hbar]]$-linearity, it extends to the linear operator $I_{\hbar}:
\diff{A}[[\hbar]]\rightarrow\diff{A}[[\hbar]]$ with the property
\begin{equation*}
I_{\hbar}=\id+\sum_{i=1}^{\infty}\hbar^n I^{(n)},
\end{equation*}
where $I^{(n)}:\diff{A}\rightarrow\diff{A}$. Therefore, the operator $I_{\hbar}$
is invertible.

Define the $\star$-product of contravariant symbols as in~\eqref{contra4},
\begin{equation*}
\hat{f}_1\hat{\star}\hat{f}_2=I_{\hbar}^{-1}(I_{\hbar}(\hat{f}_1)\star
I_{\hbar}(\hat{f}_2)),
\end{equation*}
which, of course, is equivalent to the $\star$-product of covariant symbols.

Define a formal Berezin-Toeplitz operator $F:
\diff{J}_{hol}[[\hbar]]\rightarrow\diff{J}_{hol}[[\hbar]]$ with contravariant symbol
$\hat{f}\in\diff{A}[[\hbar]]$ as follows: $F(g)=P_{\hbar}(\hat{f}g)$ for all
$g \in\diff{J}_{hol}[[\hbar]]$. We have the following statement.
\begin{proposition} The product $F=F_1F_2$ of formal Berezin-Toeplitz operators
$F_1$ and $F_2$ with contravariant symbols $\hat{f}_1$ and
$\hat{f}_2\in\diff{A}[[\hbar]]$ is a formal Berezin-Toeplitz operator with
contravariant symbol $\hat{f}= \hat{f}_1 \hat{\star}\hat{f}_2$.
\end{proposition}
\begin{proof} This is a formal analog of  lemma~\ref{Contra}.
Using Fubini theorem for formal integrals and the equality of two types
of formal integrals proved in the Appendix, we have
\begin{equation*}
P_{\hbar}(\hat{f}g)(z)=\oint_{\CC^n}e_{\hbar}(z,\bar{v})f(z,\bar{v})g(v)
e^{(\Phi(v,\bar{z})-\Phi(v,\bar{v}))/\hbar}d\mu_{\hbar}(v,\bar{v}),
\end{equation*}
where $f=I_{\hbar}(\hat{f})$. The same arguments also show that
\begin{equation*}
(F_1F_2)(g)(z)=\oint_{\CC^n}e_{\hbar}(z,\bar{v})(f_1\star f_2)(z,\bar{v})
g(v)e^{(\Phi(z,\bar{v})-\Phi(v,\bar{v}))/\hbar}d\mu_{\hbar}(v,\bar{v}),
\end{equation*}
and the statement follows.
\end{proof}
\begin{remark} Some interesting explicit formulas for the asymptotic expansion of
the Bergman kernel on $\CC$, as well as for the coefficients of the product
$\hat{\star}$, were derived in~\cite{wati}. It is also interesting that
in~\cite{wati} the deformation quantization appeared as a part of M-theory in
the string theory.
\end{remark}
\subsection*{Acknowledgments} For many inspiring discussions we are greatly indebted
to M. Flato, who passed away in November 1998. We also would like to thank
M.~Kontsevich, S.~Shatashvili and D.~Sternheimer for helpful discussions,
W.~Taylor for bringing~\cite{wati} to our attention, and A. Kogan for the
help with the figures. The work of N.R. and L.T. was partially supported by
the NSF grants DMS-97-09594 and DMS-95-00557 respectively.

\section*{Appendix}

Here we will prove the following statement.
\begin{proposition} \nonumber (Equality of two types of formal integrals)
\begin{equation*}
\int_{\CC^n}f(v,\bar{v})e^{\phi(z,\bar{z};v,\bar{v})/\hbar}d\mu_{\hbar}(v,
\bar{v})=\oint_{\CC^n}\tilde{f}_z(v,\bar{v})
e^{(\Phi(z,\bar{v})-\Phi(v,\bar{v}))/\hbar} d\mu_{\hbar}(v,\bar{v}),
\end{equation*}
where $\phi(z,\bar{z};v,\bar{v})=\Phi(z,\bar{v})+\Phi(v,\bar{z})-
\Phi(z,\bar{z})-\Phi(v,\bar{v})$ is the Calabi function of the K\"{a}hler form
and
\begin{equation*}
\tilde{f}_z(v,\bar{v})=f(v,\bar{v})e^{(\Phi(v,\bar{z})-
\Phi(z,\bar{z}))/\hbar}.
\end{equation*}
Here the formal integral in the left hand side is understood as the formal
Laplace expansion around the critical point $v=z$ of the exponential factor,
i.e., as the asymptotic of the distribution $\exp(\phi(z,\bar{z};v,\bar{v})/
\hbar)$ as $\hbar\rightarrow 0$. The formal integral in the right hand side
is understood as the asymptotic of the distribution $\exp((\Phi(z,\bar{v})-
\Phi(v,\bar{v}))/\hbar)$ as $\hbar\rightarrow 0$.
\end{proposition}
\begin{proof} We will start with case $n=1$; the precise definition of the formal
integral $\oint$ will become clear in a due course. First, assume that $f\in
C^{\infty}(\CC)_{pol}$, that the Calabi function has a single
critical point at $v=z$, and, assuming that the integral
\begin{equation*}
J(f)(\hbar;z,\bar{z})=\frac{\sqrt{-1}}{2\pi\hbar}\int_{\CC}f(v,\bar{v})
e^{\phi(z, \bar{z};v,\bar{v})/\hbar}\del^2_{v\bar{v}}\Phi(v,\bar{v})dv\wedge
d\bar{v}
\end{equation*}
exists for $\hbar>0$, consider its the asymptotic behavior
as $\hbar\rightarrow 0$. Replacement of the domain of integration by a
neighborhood $U_{\delta}=\{v~|~v\in\CC,~|v-z|<\delta\}$ of $z\in\CC$ (with
certain $\delta>0$ depending on $f$ and $\Phi$ in the support of $f$) will
result in an exponentially small error term as $\hbar\rightarrow 0$. (It
would be also sufficient to consider terms that decay faster than any
power $\hbar^N$ with $N>0$). Therefore,
\begin{eqnarray*}
J(f)(\hbar;z,\bar{z})&\cong &\lim_{\epsilon\rightarrow 0}
\frac{\sqrt{-1}}{2\pi\hbar} \int_{\epsilon<|v-z|<\delta}f(v,\bar{v})
e^{\phi(z,\bar{z};v,\bar{v})/\hbar} \del^2_{v\bar{v}}\Phi(v,\bar{v})dv\wedge
d\bar{v} \\
&\cong &\lim_{\epsilon\rightarrow 0}\frac{\sqrt{-1}}{2\pi\hbar}
\int_{\epsilon<|v-z|<\delta}\tilde{f}_z(v,\bar{v})e^{(\Phi(z,\bar{v})-
\Phi(v,\bar{v}))/\hbar}\del^2_{v\bar{v}}\Phi(v,\bar{v})dv\wedge d\bar{v},
\end{eqnarray*}
where the symbol $\cong$ stands for  equality modulo exponentially small
terms as $\hbar\rightarrow 0$.

Next, for fixed $z$, consider the following representation (cf~\cite{wati})
\begin{equation*}
e^{(\Phi(z,\bar{v})-\Phi(v,\bar{v}))/\hbar}dv\wedge d\bar{v}=d(A(\hbar; v,
\bar{v})e^{(\Phi(z,\bar{v})-\Phi(v,\bar{v}))/\hbar}d v),
\end{equation*}
where
\begin{equation*}
A(\hbar;v,\bar{v})=\sum_{n=1}^{\infty}A_n(v,\bar{v})\hbar^n,
\end{equation*}
and satisfies the equation
\begin{equation*}
-\delb_v A(\hbar;v,\bar{v})+\frac{1}{\hbar}A(\hbar;v,\bar{v})(\delb_v
\Phi(z,\bar{v}) -\delb_v \Phi(v,\bar{v}))=1,
\end{equation*}
so that
\begin{equation*}
A_1(v,\bar{v})=(\delb_v\Phi(z,\bar{v})-\delb_v\Phi(v,\bar{v}))^{-1}~
\text{and}~A_{n+1}(v,\bar{v})=A_1(v,\bar{v})\delb_v A_n(v,\bar{v}).
\end{equation*}
Applying Stokes' formula, we get
\begin{gather*}
J(f)(\hbar;z,\bar{z})\cong \frac{\sqrt{-1}}{2 \pi\hbar}
\lim_{\epsilon\rightarrow 0}\{\int_{C_{\epsilon}}\tilde{f}_z(v,\bar{v})
A(\hbar;v,\bar{v})\del^2_{v\bar{v}} \Phi(v,\bar{v})d v-
\int_{\epsilon <|v-z|<\delta} \\
A(\hbar;v,\bar{v})\del_{\bar{v}}(\tilde{f}_z\del_v\delb_v\Phi)(v,
\bar{v}))e^{(\Phi(z,\bar{v})-\Phi(v,\bar{v}))/\hbar} dv\wedge d\bar{v}\},
\end{gather*}
where $C_{\epsilon}=\{v~|~v\in\CC,~|v-z|=\epsilon\}$, and we have used that
the integral over the circle of radius $\delta$ around $z$ is exponentially
small as $\hbar\rightarrow 0$.

Repeating this procedure, we obtain
\begin{gather*}
J(f)(\hbar;z,\bar{z})\cong \frac{\sqrt{-1}}{2\pi\hbar}\lim_{\epsilon
\rightarrow 0}\sum_{n=1}^{\infty}(-1)^n\int_{C_{\epsilon}}
\diff{D}^n(\tilde{f}_z \del^2_{v\bar{v}}\Phi))(v,\bar{v})A(\hbar; v,\bar{v}) \\
e^{(\Phi(z,\bar{v})-\Phi(v,\bar{v}))/\hbar}dv,
\end{gather*}
where $\diff{D}=A(\hbar,v,\bar{v})\delb_v$. Finally, expanding $A(\hbar; v,
\bar{v})$ into power series in $\hbar$ with coefficients $A_n(v,\bar{v})$
that have singularities $a_n(z,\bar{z})(v-z)^{-n}$ as $v\rightarrow z$, and
using the generalized Cauchy formula
\begin{equation*}
\lim_{\epsilon\rightarrow 0}\frac{1}{2\pi \sqrt{-1}}\int_{C_{\epsilon}}
\frac{g(v,\bar{v})}{(v-z)^n}d v=\frac{(\del_z)^{n-1}g(z,\bar{z})}{(n-1)!},
~n=1,2,\ldots,
\end{equation*}
which is valid for any $g\in C^{\infty}(U_{\delta})$, we get an asymptotic
expansion for $J(f)(\hbar;z,\bar{z})$ in terms of a power series in $\hbar$.
Clearly, the coefficients of this expansion are given by the infinite sums,
convergent when the function $f$ and K\"{a}hler potential $\Phi$ are real-analytic.
This expansion, extended to the formal category (with respect to $\hbar$) with
coefficients given by the global sections of the bundle of infinite jets over $\CC$,
is the definition of a formal integral $\oint$. Now proposition follows by the same
arguments as in the proof of lemma~\ref{ass}.

The proof for the case $n>1$ is similar and is left to the reader.
\end{proof}

\end{document}